\documentclass[12pt,a4paper]{amsproc}
\usepackage[english]{babel}
\usepackage{latexsym}
\usepackage{amsmath}
\usepackage{amssymb}
\usepackage[all]{xy}

\usepackage{amsmath,tikz-cd}
\usepackage{amscd}
\usepackage[cp850]{inputenc}
\usepackage[mathscr]{eucal}
\newcommand{\To}{\longrightarrow}
\newcommand{\bil}[2]{{\left\langle\kern0ex #1,#2
		\kern0ex\right\rangle}}
\newcommand{\PS}{{\sf P}\hspace{-1pt}{\sf S}}

\numberwithin{equation}{section}

\newtheorem{thm}{Theorem}[section]
\newtheorem{cor}[thm]{Corollary}
\newtheorem{lem}[thm]{Lemma}
\newtheorem{prop}[thm]{Proposition}
\newtheorem{quest}{Question}[section]

\newtheorem{defin}[quest]{Definition}

\theoremstyle{remark}

\newcommand{\R}{\mathbb{R}}
\newcommand{\N}{\mathbb{N}}

\newcommand{\bA}{\mathbb A}

\newcommand{\pe}{\mathfrak p}
\newcommand{\con}{\mathfrak c}
\newcommand{\cA}{\mathcal{A}}
\newcommand{\cC}{\mathcal{C}}

\newcommand{\adef}{\begin{defin}}
\newcommand{\zdef}{\end{defin}}

\def\Ext{\operatorname{Ext}}

\def\supp{\operatorname{supp}}
\def\PB{\operatorname{PB}}
\def\PO{\operatorname{PO}}
\newcommand\restr[2]{{
		\left.\kern-\nulldelimiterspace 
		#1 
		\right|_{#2} 
}}

\def\Ext{\operatorname{Ext}}

\newcommand{\aproof}{\begin{proof}}
\newcommand{\zproof}{\end{proof}}

\title{Twisted sums of $c_0(I)$}

\author{Jes\'{u}s M. F. Castillo}
\address{Universidad de Extremadura\\Instituto de Matem\'aticas Imuex\\
Avenida de Elvas s/n\\ 06071-Badajoz\\ Spain} \email{castillo@unex.es}

\author{Alberto Salguero Alarc\'on}
\address{Universidad de Extremadura\\Instituto de Matem\'aticas Imuex\\
Avenida de Elvas s/n\\ 06071-Badajoz\\ Spain} \email{salgueroalarcon@unex.es}



\thanks{This research was supported in part by MINCIN project PID2019-103961GB-C21 and by Junta de Extremadura project IB20038.
The second-named author benefited from an FPU Grant FPU18/00990 from the Spanish Ministry of Universities.}

\begin{document}

\maketitle

\begin{abstract} The paper studies properties of twisted sums of a Banach space $X$ with $c_0(\kappa)$. We first prove a representation theorem for such twisted sums from which we will obtain, among others, the following: (a) twisted sums of $c_0(I)$ and $c_0(\kappa)$ are either subspaces of $\ell_\infty(\kappa)$ or trivial on a copy of $c_0(\kappa^+)$; (b) under the hypothesis $[\mathfrak p = \mathfrak c]$, when $K$ is either a suitable Corson compact, a separable Rosenthal compact or a scattered compact of finite height, there is a twisted sum of $C(K)$ with $c_0(\kappa)$ that is not isomorphic to a space of continuous functions; (c) all such twisted sums are Lindenstrauss spaces when $X$ is a Lindenstrauss space and $G$-spaces when $X=C(K)$ with $K$ convex, which shows tat a result of Benyamini is optimal; (d) they are isomorphically polyhedral when $X$ is a polyhedral space with property ($\star$), which solves a problem of Castillo and Papini.\end{abstract}

\section{Introduction}
We address the reader to the Preliminaries section for all unexplained terminology. Twisted sums of the form $0\To c_0\To  X \To c_0(I)\To 0$ are simultaneously simple and very complex mathematical objects. The examples of $X$ known so far are:
\begin{itemize}
\item Spaces of continuous functions $C(K_\cA)$ over the Stone compactum of a Boolean algebra generated by an uncountable almost disjoint family $\cA$ of subsets of $\N$.
\item When $|I|=\mathfrak c$, the spaces constructed in \cite{ps} (under some form of Martin's axiom), whose main feature is that they are not isomorphic to any $C(K)$-space.
\item Products and $c_0$-sums of both.
\end{itemize}

The nature and properties of such twisted sums depend on set-theoretic assumptions: in \cite{maple} it is shown that, under [{\sf CH}], there exist $2^{\aleph_1}$ non-isomorphic $C(K_\cA)$ spaces generated by families of size $\aleph_1$; while in \cite{sailing} it is shown that, under [{\sf MA} + $\aleph_1<\mathfrak c$] all such spaces are isomorphic. And in \cite{amp} it is shown that $\Ext(C(K_\cA), c_0)\neq 0$ under [{\sf CH}] while \cite{maple} proves that $\Ext(C(K_\cA), c_0)= 0$ provided $\cA$ is of size $\aleph_1$ under  [{\sf MA} + $\aleph_1<\mathfrak c$]. All this make us think that time is ripe to undertake a classification of such twisted sum spaces.

\section{Preliminaries}
\subsection{Background on exact sequences}
\par An \emph{exact sequence} of Banach spaces is a diagram of Banach spaces and (linear, continuous) operators
\begin{equation} \label{eq:seq}  \xymatrix{0 \ar[r] & Y \ar[r]^i & Z \ar[r]^q & X \ar[r] & 0} \end{equation}
so that the kernel of every operator coincides with the image of the preceding one. The middle space $Z$ is usually called a \emph{twisted sum} of $Y$ and $X$. The open mapping theorem asserts that in such a situation $Y$ is isomorphic to a closed subspace of $Z$ and $X$ is isomorphic to $Z/Y$. Even if, in general, the middle space $Z$ in (\ref{eq:seq}) where $Y$ and $X$ are both Banach spaces is only a quasi-Banach space and does not need to be a Banach space, we will not have the need for such considerations, since a deep result of Kalton and Roberts \cite{kalrob} guarantees that $Z$ is a Banach space whenever $X$ is an $\mathscr{L}_\infty$-space.

\par Two exact sequences $0\To Y\To  Z_i \To X\To 0$, $i=1,2$ are said to be \emph{equivalent} if there is an operator $u: Z_1 \to Z_2$ making commutative the diagram
\[\begin{CD} 0@>>> Y@>>> Z_1 @>>> X@>>> 0\\
	&&@| @VVuV @|\\
	0@>>> Y @>>> Z_2 @>>> X@>>> 0\\
	\end{CD}
	\]
The $3$-lemma and the open mapping theorem assure that $u$ is an isomorphism, and so this defines a true equivalence relation.  We say that the exact sequence (\ref{eq:seq}) is \emph{trivial}, or that \emph{it splits}, when it is equivalent to the exact sequence $0\To Y \To  Y \oplus X  \To X\To 0$. This happens if and only if $i(Y)$ is a complemented subspace of $Z$, and if and only if there is a continuous right-inverse for the quotient map $q$.  We shall write $\Ext(X,Y)$ for the set of equivalence classes of exact sequences (\ref{eq:seq}). Using pullbacks and pushouts (see their definitions below), this set can be endowed with a vector space structure so that the zero element is the class of trivial sequences. Therefore, $\Ext(X,Y)=0$ means that every twisted sum of $Y$ and $X$ is trivial.

\par Two basic constructions regarding exact sequences are the \emph{pullback} and the \emph{pushout}. Consider an exact sequence (\ref{eq:seq}) and two operators $T: X' \to X$ and $S: Y \to Y'$. Then there are new exact sequences forming commutative diagrams
\begin{equation} \label{eq:pbpo}
\xymatrix{0 \ar[r] & Y \ar[r]^i& Z \ar[r]^q  & X \ar[r]      & 0 &  & 0 \ar[r] & Y \ar[r]^i\ar[d]^S & Z \ar[r]^q\ar[d]^S & X \ar[r] \ar@{=}[d] & 0 \\
	0 \ar[r] & Y \ar[r]^{\overline{i}} \ar@{=}[u] & \PB \ar[r]^{\overline{q}} \ar[u]^{\overline{T}} & X \ar[r] \ar[u]^T & 0 &  & 0 \ar[r] & Y' \ar[r]^{\underline{i}}   & \PO \ar[r]^{\underline{q}}               & X \ar[r]           & 0}
\end{equation}
In the left diagram, the \emph{pullback space} is $\PB=\{(z,x')\in Z \times X': qz=Tx'\}$. The maps $\overline T$ and $\overline{q}$ are the restrictions of the canonical projections from $Z \times X'$ into $Z$ and $X'$, respectively, and $\overline{i}$ is simply the natural inclusion $\overline{i}(y) = (iy,0)$. In the right diagram, the \emph{pushout space} $\PO$ is the quotient of the direct sum $Y' \oplus Z$ by the closure of $\{(i(y), -S(y')): y\in Y\}$. The maps $\underline{S}$ and $\underline{i}$ are just the composition of the natural inclusions of $Y'$ and $Z$ into $Y' \oplus Z$, respectively,  with the quotient map, and $\underline{q}$ arises from the factorization of the operator $Y'\oplus Z \to X$, $(y',z)\mapsto q(z)$. We shall call the lower rows of the left and right diagrams in (\ref{eq:pbpo}) the \emph{pullback} and \emph{pushout} sequences, respectively. A pullback sequence splits if and only if there is a \emph{lifting} for $T$; that is, an operator $T': X' \to Z$ such that $qT'=T$. On the other hand, a pushout sequence splits if and only if $S$ \emph{extends} to $Z$, which means there is an operator $S': Z \to Y'$ so that $S'i=S$.

\subsection{$\cC$-spaces}
We will say that a Banach space is a \emph{$\cC$-space} when it is isomorphic to some space $C(K)$ of continuous functions on some compact space $K$. As usual, we identify the dual of $C(K)$ with the space $M(K)$ of Radon measures on $K$, with the variation norm. $M_1(K)$ denotes the closed unit ball of $M(K)$ endowed with the weak* topology.

\section{Everything old is new again}\label{background}
Let us begin stating the properties that every twisted sum space $0\To c_0\To  X \To c_0(I)\To 0$ must have.
\begin{prop} The twisted sum space $X$:
\begin{enumerate}
\item Is $c_0$-saturated and $c_0$-uppersaturated.
\item Is an Asplund space. Consequently, it has weak*-sequentially compact dual ball, it is weak*-extensible and has the Gelfand-Phillips property.
\item Is isomorphic to a Lindenstrauss space; consequently, it has Pe\l czy\'nski's property $(V)$.
\item $X\simeq X\oplus c_0$.
\item Is WCG if and only if is isomorphic to $c_0(I)$.
\item Is separably injective and not universally separably injective.
\end{enumerate}
\end{prop}
\begin{proof} Recall that a space is $c_0$-saturated if every closed subspace contains a copy of $c_0$; the space is $c_0$-uppersaturated \cite[Def. 2.25]{2132} if every separable subspace is contained in a copy of $c_0$ contained in the space. Both are three-space properties: $c_0$-saturation is in \cite{castgonz} while $c_0$-uppersaturation is in \cite[Prop. 6.2]{2132}. Concerning (2), the Asplund character is a consequence of (1). Now, that Asplund spaces have  weak*-sequentially
compact dual ball appears in \cite[Corollary 2]{hajo}, and the weak*-extensibility character follows from \cite{wz}. See also \cite[Proposition 6]{cgp}. This also shows that \cite[Example 0.5.6]{schl} does not require any special choice of the almost disjoint family. The Gelfand-Phillips property comes as a consequence of the weak*-sequential compactness of the dual ball, see \cite[proposition 6.8.c]{castgonz}. (3) The fact that every twisted sum of $c_0$ and a Lindenstrauss space is isomorphic to a Lindenstrauss space can be found in \cite[Theorem 5.1]{ps}, but the following argument is more informative:

\begin{lem}\label{Mideal} $c_0$ is isomorphic to an $M$-ideal in every superspace. \end{lem}
\begin{proof} Indeed, $c_0$ is an $M$-ideal in $\ell_\infty$, hence in every $\ell_\infty(I)$, in its natural position. Since $\ell_\infty(I)$ is separably automorphic \cite[\S 2.6]{2132}, it means that $c_0$ is an $M$-ideal in any $\ell_\infty(I)$ in every position. Finally, if $c_0 \subset Y\subset X$ is an $M$-ideal in $X$ then it is also an $M$-ideal in $Y$. See also \cite{hww}.\end{proof}

Assertion (3) follows immediately from that since being an $M$-ideal means that $X^*= (c_0)^* \oplus_1 (X/c_0)^*$ and thus, if $X/c_0$ is a Lindenstrauss space, then $X^* = \ell_1\oplus_1 L_1(\mu)$ and $X$ is a Lindenstrauss space. Lindenstrauss spaces have property $(V)$ and are Lindenstrauss-Pe\l czy\'nski spaces in the sense of \cite{LP}. In particular, $X$ is an $\mathcal L_\infty$-space, but recall that not all $\mathcal L_\infty$-spaces enjoy property $(V)$ or have to be a Lindenstrauss-Pe\l czy\'nski space \cite{LP}. A proof for
the more general fact that all twisted sums of $c_0(J)$ and a space with property $(V)$ have property $(V)$ can be seen in \cite{castsimo}. It is well-known that the space $c_0(I)$ is WCG and also that every copy of $c_0$ in a WCG space is complemented \cite{johnlindwcg}, which is enough to prove (4). Finally, to dispose of (5), recall that ``being separably injective" is a three-space property \cite[Prop. 1.4]{2132} and that $c_0(I)$ is separably injective, but no twisted sum of $c_0$ and $c_0(I)$ can be universally separably injective since they are never Grothendieck spaces \cite[Prop. 2.8]{2132}.
\end{proof}

The space $c_0(I)$ is the simplest twisted sum of the type we are considering. One has:

\begin{lem} $c_0(I)$ has the following properties:
\begin{enumerate}
\item Every copy of $c_0(J)$ inside $c_0(I)$ is complemented.
\item Every nonseparable subspace of density $\kappa$ contains a complemented copy of $c_0(\kappa)$.
\end{enumerate}\end{lem}
\begin{proof} The first assertion is in \cite{asg}. The second is somehow folklore and can be partially found in different places
\cite{FinolWojtowicz,orty,salinas}. Probably the deepest argument is given in \cite[Lemma 2.7]{MP}: for any subspace $X$ of $c_0(I)$ there is a decomposition $I=\cup_{i\in J} I_i$ with $|J|=|I|$ and such that if $X_i= \{x\in X: \supp x\subset I_i\}$ then $X = c_0(J, X_i)$ and the following sequences are isomorphic.
$$\xymatrix{0 \ar[r] & X \ar[r] & c_0(I) \ar[r]& c_0(I)/X \ar[r]  & 0 & \\
		0 \ar[r] & c_0(J, X_i) \ar[r] \ar@{=}[u] & c_0(J, c_0(I)) \ar[r] & c_0(J, c_0(I)/X_i) \ar[r]  & 0}
	$$
In particular it is clear that $X$ contains $c_0(J)$.\end{proof}

\subsection{Twisted sums of $c_0$ and $C(K)$ that are not $\cC$-spaces}
\par The purpose of \cite{ps} was to construct, for every Eberlein compactum $K$ of weight $\con$ a twisted sum of $c_0$ and $C(K)$ which is not a $\cC$-space. Let us recall the key points of the construction. The cardinal $\pe$ is defined in \cite[11 D]{fremlin}.
\begin{itemize}
	\item First, Theorem 3.4 in \cite{ps} produces, under $\pe=\con$, a very special almost disjoint family $\cA$ of subsets of natural numbers and size $\con$. Its existence relies on the fact that, since $K$ is Eberlein, $M_1(K)$ is a sequentially compact space and $|M_1(K)|=\con$.
	\item Next, one needs to find a subspace of the form $\Sigma = \{ \sigma_\xi, -\sigma_\xi: \xi<\con\} \cup \{0\}$ inside $M_1(K)$ which is homeomorphic to $\bA(\con)$, the one point compactification of a discrete space of cardinality $\con$. This is done in \cite[Lemma 4.1]{ps}.
	\item Then, the authors consider a compactum $L$ defined as the adjunction space of $K_\cA$ and $M_1(K)$ under an appropriate map $\psi: K_\cA' \to \Sigma \subseteq M_1(K)$. Note that $L$ is a \emph{countable discrete extension} of $M_1(K)$; namely, it contains $M_1(K)$ so that $L \setminus M_1(K)$ is a copy of a discrete countable space. The desired twisted sum $Z(K,\Sigma)$ appears as the pullback sequence of the diagram below:
		\[ \xymatrix{0 \ar[r] & c_0 \ar[r]^{e} & C(L) \ar[r]^{r}& C(M_1(K)) \ar[r]  & 0 & \\
				0 \ar[r] & c_0 \ar[r] \ar@{=}[u] & Z(K, \Sigma)\ar[r] \ar[u] & C(K) \ar[r] \ar[u]_{\phi} & 0}
		\]
	Here $e$ stands for the natural extension operator letting $e(f)(\mu)=0$ for every $\mu\in M_1(K)$, $r$ is the restriction map and $\phi: C(K) \to C(M_1(K))$ is the canonical evaluation map $\phi(f)(\mu) = \bil{\mu}{f}$.
\end{itemize}

\par Having in mind the previous requisites, one can substantially enlarge the list of compacta $K$ for which there exists a twisted sum of $c_0$ and $C(K)$ that is not a $\cC$-space.

\begin{thm} $[\pe =\con]$ There is a twisted sum $0 \To c_0 \To Z(K, \Sigma) \To C(K) \To 0$ in which $Z(K, \Sigma)$ is not a $\cC$-space provided $K$ has weight $\con$ and belongs to any of the following classes:
	\begin{enumerate}
		\item Corson compacta with property (M).
		\item Separable Rosenthal compacta.
		\item Scattered compacta of finite height.
	\end{enumerate}
\end{thm}
\begin{proof} Recall that a compactum $K$ has property (M) if every regular (Radon) measure on $K$ has separable support. A good place to look for the standard facts about Corson compacta is \cite[Ch. 6]{negrepontis}.  In particular, we mention here that separable Corson compacta are metrizable, and that Corson compacta are sequentially compact. Moreover, for a Corson compacta $K$, having property (M) is equivalent to $M_1(K)$ being Corson.
Under [{\sf MA} + $\aleph_1<\mathfrak c$], every Corson compactum has property (M), but this is not true under [{\sf CH}]. As a consequence, it is readily seen that $M_1(K)$ is sequentially compact, and also that $|M_1(K)|=\con$ because $K$ has $\con$ many metrizable subspaces and each of them can only support $\con$ many measures. Now, in order to find a copy of $\bA(\con)$ inside $M_1(K)$ when $K$ is Corson, we note that  \cite[Lemma 4.1]{ps} is also valid for Corson compacta.

\par Now we assume that $K$ is a separable Rosenthal compactum. The properties of Rosenthal compacta needed to carry out the proof are contained in the survey paper by Marciszewski \cite{marciszewski}. Especially, we infer that $M_1(K)$ is sequentially compact and $|M_1(K)|=\con$ from the fact that $M_1(K)$ is separable Rosenthal provided $K$ is. Finally, $M_1(K)$ contains a copy of $\bA(\con)$ by virtue of a theorem of Todor\v cevi\'c \cite[Th.9]{todorcevic}.
\par Last, we assume that $K$ is a scattered compacta of finite height and weight $\con$. Then $C(K)$ is Asplund, and so $M_1(K)$ is sequentially compact. Also, $M(K)$ is isomorphic to $\ell_1(\con)$, so $|M_1(K)|=\con$. Finally, the assumption of finite height enables us to find a copy of $\bA(\con)$ inside $K$, and therefore inside $M_1(K)$, thanks to \cite[Lemma 6.4]{amp}. We point out that $\con$ is a regular cardinal under $\pe=\con$ because $\pe$ is regular \cite[21K]{fremlin}.
\end{proof}

\par It remains open whether $\beta\omega$ can or cannot be added to the list. However, $\Ext(\ell_\infty, c_0)\neq 0$, as shown in \cite[\S 2.1]{CC04}. Another question that stems from \cite{ps} is whether the spaces $Z(K, \Sigma)$ are subspaces of $\ell_\infty$. This is of course evident if $K$ is separable, since in that case $C(K)$ is a subspace of $\ell_\infty$ and ``being a subspace of $\ell_\infty$'' is a three-space property. The following theorem settles this matter for Corson compacta.
\begin{thm} $[\pe=\con]$ For every Corson compacta $K$ of weight $\con$ and having property (M) there is a twisted sum
$\xymatrix{0 \ar[r] & c_0 \ar[r] & Z(K,\Sigma) \ar[r] & C(K) \ar[r] & 0}$ in which $Z(K,\Sigma)$ is a subspace of $\ell_\infty$ and not a $\cC$-space.
\end{thm}
\begin{proof} The hypothesis on $K$ means, once again, that $M_1(K)$ is Corson, and then the main result of \cite{vwz} implies that $C(K)$  possesses a Marku\v sevi\v c basis $\{(f_\xi, \mu_\xi) \in C(K) \times M(K): \xi<\con\}$; that is, a bounded biorthogonal system so that $\{f_\xi: \xi<\con\}$ separates the points in $M(K)$ and $\{\mu_\xi: \xi<\con\}$ separates the points in $C(K)$. We can further assume that $\|f_\xi\|\leq 1$ and $\|\mu_\xi\|\leq 1$ for every $\xi<\con$. Taking into account the previous considerations, it is clear that $\Sigma = \{\mu_\xi, -\mu_\xi: \xi<\con\}\cup\{0\}$ is a copy of $\bA(\con)$ inside $M_1(K)$. Therefore, when we produce the countable discrete extension $L$ of $M_1(K)$ that origines the space $Z(K, \Sigma)$ we obtain that the set $\{\mu_\xi, -\mu_\xi: \xi<\con\} \cup \{0\}$ is contained in the closure of the countable discrete set $L \setminus M_1(K)$, which we now identify with $\N$. As a consequence, the sequence of functionals $(\delta_n)_{n=1}^\infty$, where $\delta_n(f,g)= f(n)$ is norming, since if $(f,g)\in Z(K, \Sigma)$ satisfies $\delta_n(f,g)=0$, then $f(\mu_\xi)=\bil{\mu_\xi}{g} = 0$, from where it follows that $g=0$.
\end{proof}

\par It is time to return to our business and reexamine the case $K=c_0(I)$ for $|I|=\con$. In the light of our previous theorem, it is clear that choosing $\Sigma = \{ 1_i, -1_i: i\in I\} \cup \{0\}$ in the unit ball of $\ell_1(\con)$ produces a twisted sum $Z(K, \Sigma)$ of $c_0$ and $c_0(\con)$ which is a subspace of $\ell_\infty$ and not a $\cC$-space. Such twisted sum space will be denoted as $\PS$ in the sequel. We point out the fact that $\PS \oplus c_0(\con)$ is not a $\cC$-space either: this can be seen by considering $J \subset I$ so that $|J|=|I\setminus J|=\con$ and showing that, when  $\Sigma = \{ 1_j, -1_j: j\in J\} \cup \{0\}$, the space $Z(K, \Sigma)$ is isomorphic to $\PS \oplus c_0(\con)$.

\section{New skin for the old ceremony}
Let us fix an infinite cardinal number $\kappa$ and a compactum $K$. A \emph{$\kappa$-discrete extension} of $K$, or $\kappa$-DE for short,  is another compactum $K\cup \kappa$ which contains a copy of $K$ and a disjoint copy of a discrete set of size $\kappa$. It is clear that if $L=K\cup \kappa$ is a $\kappa$-DE, then there is an exact sequence
$\xymatrix{0 \ar[r] & c_0(\kappa) \ar[r]^e & C(L) \ar[r]^r & C(K) \ar[r] & 0}$
 where $e$ is the plain extension operator such that $e(f)(k)= 0$ for all $k\in K$, and $r$ is the plain restriction map.
\begin{thm} \label{th:functor} For every exact sequence $\xymatrix{0 \ar[r] & c_0(\kappa) \ar[r]^j & Z\ar[r]^q  & X \ar[r] & 0}$ there is a $\kappa$-discrete extension $L$ of $B_{X^*}$ and a commutative diagram:
	\begin{equation}\label{eq:main} \xymatrix{0 \ar[r] & c_0(\kappa) \ar[r]^e & C(L) \ar[r]^r& C(B_{X^*}) \ar[r]  & 0 & \\
0 \ar[r] & c_0(\kappa) \ar[r]^j \ar@{=}[u] & Z\ar[r]^q \ar[u]^u & X \ar[r] \ar[u]_{\phi_X} & 0}
			\end{equation}
where $\phi_X: X \to C(B_{X^*})$ is the canonical evaluation map $\phi_X(x)(x^*)=\bil{x^*}{x}$.
\end{thm}
\begin{proof} The dual of the original sequence, $0 \To X^* \To Z^* \To \ell_1(\kappa) \To 0$ obviously splits. Therefore, for every $\alpha\in\kappa$ there is $z^*_\alpha \in Z^*$ such that $j^*z^*_\alpha = e^*_\alpha$, where $\{e^*_\alpha:\alpha\in \kappa\}$ denotes the canonical basis of $\ell_1(\kappa)$. This implies that any weak*-cluster point of $\{z^*_\alpha\}_{\alpha\in \kappa}$ belongs to $j[c_0(\kappa)]^\perp = q^*(X^*)$,  since $\bil{z_\alpha^*}{j x} = \bil{e_\alpha^*}{x} \to 0$ for every $x\in c_0(\kappa)$. On the other hand, it is clear that no $z_\alpha^*$ belongs to $j[c_0(\kappa)]$. Consider $\varepsilon >0$ so that $\|z^*_\alpha\|>\varepsilon$ and set $L = q^*[\varepsilon B_{X^*}] \cup \{z_\alpha^*: \alpha\in \kappa\}$ endowed with the weak* topology of $Z^*$. Since $q^*[\varepsilon  B_{X^*}]$ is homeomorphic to $B_{X^*}$, $L$ can be readily identified with a  $\kappa$-discrete extension of $B_{X^*}$, which we call again $L$. The map $u: Z \to C(L)$ defined by $u(z)(x^*)=\bil{q^*x^*}{z}$ and $u(z)(z_\alpha^*)=\bil{z_\alpha^*}{z}$ makes diagram (\ref{eq:main}) commutative.\end{proof}

\begin{defin} Let $X$ be a Banach space. A discrete extension $L$ of $B_{X^*}$ is \emph{realizable inside $X^*$} if the canonical embedding $B_{X^*} \hookrightarrow X^*$ can be extended to an embedding $L \hookrightarrow X^*$.
\end{defin}

The standard splitting criterion can be reformulated now in the language of discrete extensions:

\begin{prop} The lower sequence in (\ref{eq:main}) is trivial if and only if $L$ is realizable inside $X^*$.
\end{prop}\begin{proof} It is a consequence of the fact that the functors $\bigcirc^*: \text{Ban}_1 \rightsquigarrow \text{Comp}$ and $C(\cdot): \text{Comp} \rightsquigarrow \text{Ban}_1$ are mutually adjoint; i.e., there is a natural identification $\mathcal C(L, B_{X^*}) = \mathfrak L(X, C(L))$. \end{proof}

\section{Both sides now}

The representation in Theorem \ref{th:functor} can be used to solve several problems in the area:

\subsection{Twisted sums of $c_0(\kappa)$ and Lindenstrauss spaces are isomorphically Lindenstrauss}

The proof for $\kappa=\aleph_0$ has been already obtained as (3) in Section \ref{background}. We present now a different more general argument:

\begin{thm}\label{isomLind} Every twisted sum of $c_0(\kappa)$ and a Lindenstrauss space is isomorphically Lindenstrauss.
\end{thm}
\begin{proof} Consider $Z$ a twisted sum of $c_0(\kappa)$ and $X$ and let $L$ be a $\kappa$-DE extension of $B_{X^*}$ generating a diagram (\ref{eq:main}). We shall write $\phi$ instead of $\phi_X$ for the evaluation map $X \to C(B_{X^*})$. Given $\mu\in M(B_{X^*})$ and $x\in X$, $\int_{B_{X^*}} \phi(x) \, d\mu = \phi(x) (\phi^*\mu)$, which justifies calling $\phi^*\mu$ the \emph{barycenter} of $\mu$. Complete diagram (\ref{eq:main}) to get

$$\xymatrix{
&C(L)/Z \ar@{=}[r]& C(B_{X^*})/X\\
c_0(\kappa) \ar[r]^e & C(L)\ar[u] \ar[r]^r& C(B_{X^*}) \ar[u]\\
c_0(\kappa) \ar[r] \ar@{=}[u] & Z\ar[r] \ar[u]^u & X  \ar[u]_{\phi} }$$
where we have omitted the initial and final $0$ in the exact sequences above. Consider the dual sequence $\xymatrix{	0 \ar[r] & Z^\perp \ar[r] & M(L) \ar[r]^{u^*} & Z^* \ar[r] & 0}$
of the middle column, where $u^*(\mu)=\restr\mu{\kappa} + \phi^*\mu$. Every element $[\mu] \in Z^*$ is determined by $\restr\mu\kappa$ together with the barycenter of $\restr{\mu}{B_{X^*}}$; and its norm can be computed as
	\[ \|[\mu]\| = \|\restr\mu\kappa\| + \inf_{\nu \in \ker u^*} \|\restr{(\mu-\nu)}{B_X^*}\|
	\]
One clearly has $\|\restr{(\mu-\nu)}{B_X^*}\| \geq \|\phi^*(\restr\mu{B_{X^*}})\|$ whenever $\nu \in \ker u^*$. To see the reverse inequality, extend homogeneously the embedding $x^* \mapsto \delta_{x^*}$ of $B_{X^*}$ into $M(B_{X^*})$ to produce a map $s: X^* \to M(B_{X^*})$ which is a selection for $\phi^*$. Therefore, the choice $ \nu = \restr{\mu}{B_{X^*}} - s(\phi^*\mu)$ yields $\|\restr{(\mu-\nu)}{B_X^*}\|  = \|\phi^*(\restr{\mu}{B_X^*})\|$. As a consequence, $Z^*$ is isometrically isomorphic to $\ell_1(\kappa) \oplus_1 X^*$, hence $Z(L)$ is a Lindenstrauss space provided $X$ is.
\end{proof}

\subsection{Twisted sums that are $G$-spaces}

There are other types of Lindenstrauss spaces in the literature; among them the $M$-spaces and $G$-spaces. The former are the sublattices of $C(K)$-spaces while the latter are the sublattices $X$ of $C(K)$-spaces so that for a certain set of triples $A= \{(k_\alpha^1, k_\alpha^2,\lambda_\alpha)\in K\times K\times \R \}$ one has $X=\{f\in C(K): f(k_\alpha^1)=\lambda_\alpha f(k_\alpha^2)\}$ for all $\alpha \in A$. We can improve Theorem \ref{isomLind} in some cases:

\begin{prop} Let $K$ be a convex compact. Every twisted sum of $c_0(\kappa)$ and $C(K)$ is a $G$-space.\end{prop}
\begin{proof} Diagram (\ref{eq:main}), which in this case is
$$\xymatrix{0 \ar[r] & c_0(\kappa)  \ar[r]& C(M_1(K) \cup \kappa ) \ar[r]& C(M_1(K)) \ar[r]& 0 \\
		0 \ar[r] & c_0(\kappa) \ar[r] \ar@{=}[u] & Z \ar[r] \ar[u] & C(K) \ar[r] \ar[u]_{\phi} & 0}$$
yields
\begin{eqnarray*} Z&\simeq  &\left \{ (f, g) \in C \left (M_1(K) \cup \kappa \right )\oplus C(K) :  f(\mu) = \mu (g) \quad \forall \mu \in M_1(K)\right \}\end{eqnarray*}

Let $\delta: C(K)\To C(M_1(K))$ be the canonical embedding. Since $K$ is convex, every $\mu\in M_1(K)$ admits a barycenter $m(\mu)\in K$, namely, a point such that $\bil{\mu}{g}= g(m(\mu))$ for every $g\in C(K)$ and $\mu\in M_1(K)$ (see \cite[Proposition 1.1]{phelps}). Let $\biguplus$ denote disjoint union so that one has $C \left (M_1(K) \cup \kappa \right )\oplus C(K) = C\left((M_1(K) \cup \kappa \right ) \biguplus K)$ with correspondence
$(h_1, h_2) \longleftrightarrow h$. Thus,
\begin{eqnarray*} Z&\simeq  &\left \{ h \in C\left((M_1(K) \cup \kappa \right ) \biguplus K) :  h_1(\mu) = h_2(m(\mu)) \quad \forall \mu \in M_1(K)\right \}\end{eqnarray*}
which is the $G$-space generated by the set of triples $\{(\mu, m(\mu), 1): \mu \in M_1(K)\}$.\end{proof}

Benyamini proved in \cite{benyG} that separable $G$-spaces are isomorphic to $\mathcal C$-spaces. Benyamini's result is optimal, at least consistently, because

\begin{cor} Under $\mathfrak p = \mathfrak c$, there exists nonseparable $G$-spaces not isomorphic to any $\mathcal C$-space.\end{cor}
\begin{proof} Indeed, let $K$ be a nonmetrizable Eberlein compact. As can be seen for instance in \cite[Theorem 1.2.4]{fab}, $M_1(K)$ is a convex  Eberlein compact. Thus, by \cite{ps}  there is an exact sequence $\xymatrix{0\ar[r]& c_0 \ar[r]&\PS(M_1(K)) \ar[r]& C(M_1(K)) \ar[r]&0}$ in which $\PS(M_1(K))$, which is a $G$-space, is not however isomorphic to a $\mathcal C$ space.\end{proof}


\subsection{Twisted sums of $c_0(I)$ are isomorphically polyhedral}

A Banach space is polyhedral if every finite dimensional (equivalently, $2$-dimensional) subspace is isometric to a subspace of $\ell_\infty$. A Banach space is isomorphically polyhedral if it is isomorphic to a polyhedral space. For instance, $c_0$ is polyhedral and $c$ is isomorphically polyhedral but not polyhedral; in fact, no $C(K)$-space can be polyhedra in its natural $\sup$ norm. Many times and in many places the first author asked whether ``being isomorphically polyhedral" is a three-space property. We now provide a solution to this problem.
Let us first recall one of the most effective criteria for obtaining isomorphical polyhedrality. Given a Banach space $X$, a subset $B$ of $S_{X^*}$ is called a \emph{boundary} if for every $x\in X$ there is $x^*\in B$ so that $\bil{x^*}{x}=\|x\|$. A boundary $B$ of $X$ is said to have \emph{property} $(\star)$ if for every weak*-cluster point $x^*$ of $B$ we have $x^*(x)<1$ whenever $\|x\|=1$.

\begin{prop} \emph{\cite[Prop. 5]{fpst}} Every Banach space admitting a boundary with property $(\star)$ is polyhedral.
\end{prop}

It is an open question whether if every polyhedral space possesses a boundary with property $(\star)$. A discussion on the topic, together with some partial results in this direction, can be found in \cite{troy}. We have:

\begin{thm} Let $X$ be a Banach space admitting a boundary with property $(\star)$. Every twisted sum of $c_0(\kappa)$ and $X$ admits, under an equivalent renorming, a boundary with property $(\star)$. In particular, it is isomorphically polyhedral.
\end{thm}
\begin{proof} Let $\mathcal B$ be a boundary for $X$ with property $(\star)$, and consider $Z$ a twisted sum of $c_0(\kappa)$ and $X$ renormed, as in (\ref{eq:main}), as a subspace of $C(B_{X^*}\cup \kappa)\oplus_\infty X$. We show that
	\[ \widehat{\mathcal B }= \{\delta_\alpha: \alpha < \kappa\} \cup \mathcal B \subseteq Z^* \]
is a boundary for $Z$ having property $(\star)$. Indeed, if $z^*$ is a  weak*-cluster point of  $\widehat{\mathcal B }$ then it vanishes on $c_0(\kappa)$. We can therefore consider $z^*$ as a weak*-cluster point of $\mathcal B$ in $X^*$, which automatically implies $z^*(f,x)=\bil{z^*}{x}<1$ for every point $(f,x)\in Z$ of norm 1. \end{proof}

\begin{cor} Every twisted sum of $c_0(\kappa)$ and $c_0(I)$ is isomorphically polyhedral.
\end{cor}

This solves a problem repeatedly posed by Castillo and Papini \cite{troy,castpapi}.

\subsection{A dichotomy for twisted sums of $c_0(I)$}

Let $L = B_{\ell_1(I)} \cup \kappa$. We define
$\kappa_\perp = \{f\in Z(L):f(\alpha)=0 \; \forall \alpha \in \kappa\}$ and $\beta = \mathrm{dens} (\kappa_\perp) $. One has:

\begin{thm} Every twisted sum of $c_0(\kappa)$ and $c_0(I)$ is either a subspace of $\ell_\infty(\kappa)$ or it is trivial on a copy of $c_0(\kappa^+)$. \end{thm}
\begin{proof} The restriction map $R: \kappa_\perp\To C(B_{\ell_1(I)})$ is an isomorphism, and thus $R[\kappa_\perp]$ is a subspace of $c_0(I)$. Let $(f_i)_{i\in \beta}$ be a separating family of functionals on $\kappa_\perp$. One has $\kappa\leq \beta \leq |I|$.
The set $\{f_i, \delta_j\}_{i\in \beta, \, j\in \kappa}$ is separating in $Z$, which implies $Z$ is a subspace of $\ell_\infty(\max\{\beta, \kappa\})$. If $\beta =\kappa$ then $Z$ is a subspace of $\ell_\infty(\kappa)$. If $\beta>\kappa$ then $R[\kappa_\perp]$ is a subspace of $c_0(I)$ of density $\beta$ and therefore contains a copy of $c_0(\kappa^+)$. Since that copy is complemented in $c_0(I)$ by the result of Granero \cite{asg}, it is complemented in the twisted sum space too.\end{proof}

\section{Yesterday, once more}

The paper \cite{ps2} contains a consistent negative solution to the longstanding open problem about complemented subspaces of $C(K)$-spaces.
If the space $\PS$ were complemented in some $C(K)$-space it would also be a solution. The paper \cite{ps2} contains an additional interesting example (cf. theorem 6.2): an exact sequence $0\To c_0\To C(\PS_2)\To Q\To 0$ in which $Q$ is not isomorphic to a $C(K)$-space. This is rather interesting in combination with Lemma \ref{Mideal}. Observe that if $J$ is an $M$-ideal of $C(S)$ then $C(S)/J$ is a $C(K)$-space. Thus, on  one hand one gets that that result is no longer valid for \emph{isomorphic} copies and on the other that $c_0$ is isomorphic to an $M$-ideal of $C(\PS_2)$ but it is not itself an $M$-ideal. Consider the pullback diagram
$$\xymatrix{
0\ar[r]&c_0 \ar[r]& C(\Delta_{\mathcal M}) \ar[r]& c_0(I)\ar[r]&0\\
0\ar[r]&c_0 \ar@{=}[u]\ar[r]& \PB\ar[u] \ar[r]& \PS\ar[u] \ar[r]&0\\
&& c_0 \ar@{=}[r]\ar[u]& c_0 \ar[u]&\\
}$$

Now, either $\PB$ is a $C(K)$-space or it is not.
\begin{itemize}
\item If $\PB$ is a $C(K)$-space then the middle horizontal sequence provides another example of $c_0$ inside some $C(K)$-space such that the quotient $C(K)/c_0 = \PS$ is not a $C(K)$-space
\item If $\PB$ is not a $C(K)$-space then the middle vertical sequence provides a twisted sum of $c_0$ and $C(\Delta_{\mathcal M})$ that is not isomorphic to a $C(K)$-space. Moreover, the diagonal pullback sequence $0\To c_0\To \PB \To c_0(\Gamma)\To 0$ would be another example of twisted sum of $c_0$ and $c_0(I)$ not isomorphic to a $C(K)$-space
\end{itemize}

\end{document}